\documentclass[pagenumbers]{stabs2021_plain}
\usepackage{dsfont}
\usepackage{flushend}
\usepackage{comment}

\usepackage{xcolor}
\definecolor{OliveGreen}{cmyk}{0.64,0,0.95,0.40}

\newcommand{\R}{\mathbb{R}}
\newcommand{\N}{\mathbb{N}}
\newcommand{\eps}{\varepsilon}

\title{A new stochastic framework for ship capsizing}

\author{
  Manuela L.~Bujorianu,\\\textit{Maritime Safety Research Centre, University of Strathclyde, UK},
  \href{mailto:luminita.bujorianu@strath.ac.uk}{Luminita.Bujorianu@strath.ac.uk}
  \and
  Robert S.~MacKay,\\\textit{Mathematics Institute,
        University of Warwick, UK},
  \href{mailto:r.s.mackay@warwick.ac.uk}{R.S.MacKay@warwick.ac.uk}
  \and
  Tobias Grafke,\\\textit{Mathematics Institute,
        University of Warwick, UK},
  \href{mailto:T.Grafke@warwick.ac.uk}{T.Grafke@warwick.ac.uk}
  \and
  Shibabrat Naik,\\\textit{School of Mathematics, University of Bristol, UK},
  \href{mailto:S.Naik@bristol.ac.uk}{S.Naik@bristol.ac.uk}
  \and  
  Evangelos Boulougouris,\\\textit{Maritime Safety Research Centre, University of Strathclyde, UK},\\
  \href{mailto:Evangelos.Boulougouris@strath.ac.uk}{Evangelos.Boulougouris@strath.ac.uk}
}

\begin{stabsabstract}
We present a new stochastic framework for studying ship capsize.  It is a synthesis of two strands of transition state theory.  The first is an extension of deterministic transition state theory to dissipative non-autonomous systems, together with a probability distribution over the forcing functions.  The second is stochastic reachability and large deviation theory for transition paths in Markovian systems.  In future work we aim to bring these together to make a tool for predicting capsize rate in different stochastic sea states, suggesting control strategies and improving designs.
\end{stabsabstract}

\keywords{Transition State Theory, Transition Path Theory, Flux-over-saddle, Markov Models, Stochastic Reachability, Capsize Probability, Large Deviations.}

\begin{document}

\makestabstitle

\section{Introduction}
A new stochastic framework for studying ship capsize is presented for a general class of sea states exceeding mere regular waves. It has two strands, both starting from transition state theory~\cite{truhlar-garrett-klippenstein:1996, waalkens_wigners:2008} (in which we include transition path theory~\cite{vanden-eijnden:2006}). The common outcomes are survivability probabilities \cite{surviv_Prob:2010}, 
the probability rate for capsize and the most likely paths to capsize.

First is a formulation of capsize for given forces and moments as functions of time and state, leading from given initial condition to a deterministic time to capsize (infinite if no capsize) and hence from a probability distribution on initial conditions to a distribution of times to capsize.  This is based on a proposed extension of the ``flux over a saddle'' paradigm~\cite{mackay:1990} to include dissipation and non-autonomous forcing.
To take into account uncertainty about the forcing, we consider probability distributions over forcing functions (together with initial conditions) and aim to deduce the survivability probability, the probability rate for capsize per unit time as a function of time, and the most likely paths to capsize.

The second is stochastic reachability theory~\cite{bujorianu:2012} and large deviation theory~\cite{dembo-zeitouni:2010} for transitions of Markovian processes in continuous state-space. Defining unsafe regions to be avoided in state space, we can formally write down the probability of observing trajectories that start at normal conditions and reach an unsafe set. We compute this probability asymptotically in various limits via large deviation theory to avoid inefficient sampling problems. This allows us to efficiently explore stochastic capsize events and obtain the probability rate for capsize per unit time (the reachability from stochastic reachability), the survivability probability (the viability from stochastic reachability) and the most likely path to capsize (the large deviation minimiser).

An important strand that we do not address here is how to pass from a given incident field of wave, wind and current to the resulting forces and moment on the ship.
Another is how to formulate safety conditions for the operation or design of a ship, because that would depend on the above mapping.
Another is the formulation of control strategies for a ship's captain to avoid capsize, such as change of speed or course.

We contrast our framework with previous approaches, represented for example by~\cite{umeda_model:1995,thompson:1997,mccue_applications:2011,neves_contemporary:2011,fossen_parametric:2011}.  One is the study of response to periodic forcing, including the resulting bifurcations between attractors, e.g.~\cite{spyrou_dynamic:1996}; this gives very useful insights but real-world forcing is not periodic.  
 Periodic forcing has been combined with white noise 
 ~\cite{lin_chaotic:1995, falzarano_stochastic:2011} but this is still a limited perspective.  
 Rough seas are typically modelled as a train of random waves from some probability distribution and hence capsize in rough sea requires inclusion of more general stochastic processes~\cite{perez:2006}.
 Statistical approaches include extreme value theory \cite{leadbetter-lindgren-rootzen:2012, belenky_extreme:2019}, where a universal form is derived for extreme values from various types of stochastic process, but the known results require quite strong hypotheses and the approach to the asymptotic regime can be very slow. 

Here is the structure of the paper.  We begin with a rapid statement about our ship models in section~\ref{sec:ship-models}.  In section~\ref{sec:flux-over-saddle} we explain the flux over a saddle paradigm and its adaptation here.  Then we summarise the use of stochastic reachability theory and large deviation theory in section~\ref{sec:trans-state-theory}.  We bring these two strands together into a synthesis in section~\ref{sec:interconnections} and conclude in section~\ref{sec:conclusion}.

\section{Ship models}
\label{sec:ship-models}

\begin{figure}[htb]
    \centering
    \includegraphics[width=\columnwidth]{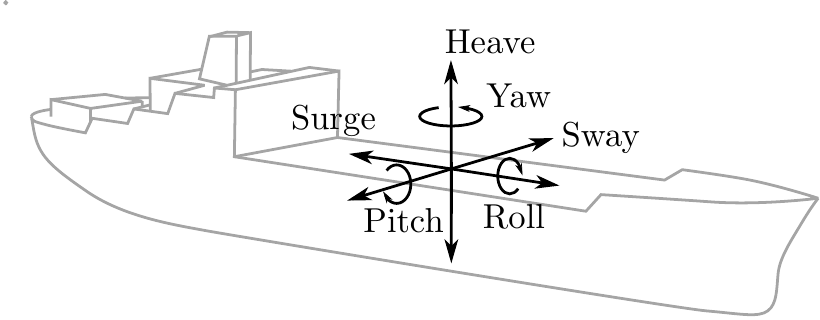}
    \caption{Degrees of freedom for ship motion.}
    \label{fig:ship-dofs}
\end{figure}
Following standard practice~\cite{belenky-spyrou-walree-etal:2019,  Evangelos:2020}, we consider a ship as a rigid body with six degrees of freedom:~roll, pitch, yaw, heave, surge and sway, subject to external forces and moments, as sketched in Fig.~\ref{fig:ship-dofs}.  Each degree of freedom consists of a configuration variable and a velocity or momentum.  The ship has an associated $6\times 6$ inertia matrix, giving the kinetic energy as a function of the state of the ship (including added mass effects for the surrounding fluid).
In addition, we use phenomenological damping forces and moments.  The result is a coupled system of 6 second-order differential equations, or equivalently of 12 first-order equations.

\section{Flux over a saddle}
\label{sec:flux-over-saddle}

Just as a continental divide separates points from which water flows to different oceans, and it consists in a set of points whose gradient trajectories flow to a saddle, the set of points whose trajectories flow to a saddle plays a key role in understanding capsize.

The starting point for the ``flux over a saddle'' paradigm is an autonomous Hamiltonian system with a saddle point having just one downhill dimension for the energy function \cite{mackay:1990}.  The saddle then possesses a ``centre manifold'' of codimension-2 in the state space ({\em codimension-2} means it has 2 dimensions less than the total state space), representing the set of initial conditions whose trajectories remain close to the saddle.  The centre manifold has a forwards contracting manifold $W^+$ of codimension-1 (commonly called its ``stable manifold'') representing states whose forward trajectories converge to trajectories on the centre manifold.  It also has a backwards contracting manifold $W^-$ (``unstable manifold'') representing states whose backwards trajectories converge to trajectories on the centre manifold.
Also the centre manifold can be spanned by a codimension-1 manifold, separated into two parts by the centre manifold (in the same way that the equator can be spanned by the surface of the earth, separated into two hemispheres).  This manifold divides the state space into two parts, corresponding to the two sides of the saddle.  To get from one side to the other, a trajectory has to cross it.  The two parts correspond to the two directions of crossing.  There is some arbitrariness in the choice of the dividing manifold, but it makes only a minor difference to when a trajectory is declared to have crossed.  The manifold $W^+$ separates the region that will make the transition from the region that will not.  Thus to find the region that will capsize one has to follow $W^+$ backwards in time.  If it avoids a core around the upright state of the ship, then the ship can be considered safe from capsize.
Use of the flux over a saddle picture in the context of ship capsize was suggested by \cite{naik-ross:2017}.

\begin{figure}[htb]
   \centering
   \includegraphics[width=\columnwidth]{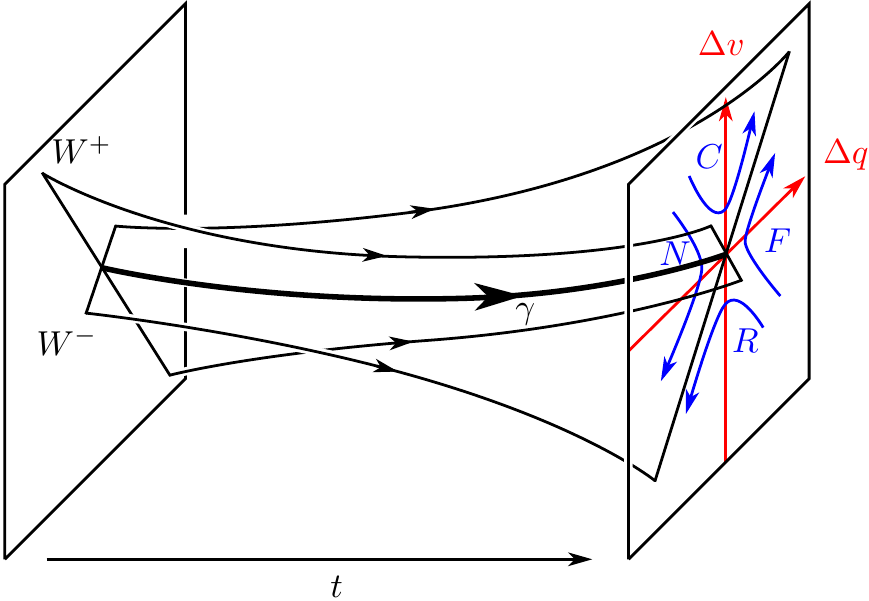} 
   \caption{Schematic of the geometry in extended state space:~the curve $\gamma$ represents a bundle of trajectories of dimension 11 that remain at the capsize threshold; deviations from it are denoted $\Delta q$ in position/attitude and $\Delta v$ in velocity; $\gamma$ has codimension-one forwards and backwards contracting submanifolds $W^\pm$, dividing the space into four sectors, labelled $C$ for capsize, $N$ for near-capsize, $R$ for righting, and $F$ for failed righting.  
   }
   \label{fig:sketch}
\end{figure}

Although originally developed in the context of Hamiltonian systems (for transition state theory of chemical reactions), the above picture persists for weak dissipation.  Furthermore, it generalises from systems with a saddle to ones with what we call a ``saddle manifold'', being a ``normally hyperbolic'' submanifold of codimension-2 with one forwards and one backwards contracting dimension ({\em normally hyperbolic} means all tangential contraction in either direction of time is slower than normal contraction in that direction of time).  We believe this is the case for a large range of realistic parameters for the standard ship models introduced in section~\ref{sec:ship-models}, with the saddle manifold being specified roughly as zero roll-velocity and a critical roll-angle as a function of all the other variables and their velocities (actually, two saddle manifolds, for port and starboard roll, and the interaction of their contracting manifolds is important).  Lastly, the framework has a version for non-autonomous systems, as is needed for periodic or more general time-dependent forcing. To describe this, we extend the 12-dimensional state space by adding time as a 13th variable. Then, if the time-dependence is not too strong, the centre manifold of the saddle has a locally unique continuation as a normally hyperbolic manifold of dimension 11 in the extended state space, that we denote by $\gamma$.  Its backwards and forwards contracting manifolds persist too, denoted by $W^\pm$.  The dividing manifold can be continued to a dividing manifold in the extended state-space.  Hence capsize for a dissipative, non-autonomous system is described by passage over this generalised saddle.  A 3D sketch of the situation is given in Figure~\ref{fig:sketch}.

We define the time $T$ to capsize to be the time until the first intersection with the dividing manifold, with the convention that $T=\infty$ if it is never reached.  Thus from a probability distribution over initial states, we obtain a probability distribution for the time $T$ to capsize.  Its derivative is the probability rate for capsize at time $T$.

In addition to probability distribution over initial conditions, we are interested in taking probability distributions over the forcing functions.  Then we want to compute features of the probability distribution of the time $T$ to capsize, in particular what is its rate as a function of $T$, and what is the probability of eventual capsize?  More broadly, what are the most likely routes to capsize?  How do all these depend on the probability distributions for the forcing functions and initial conditions, and on the parameters of the ship model? 

\section{Stochastic Reachability and Large Deviations}
\label{sec:trans-state-theory}

\begin{figure*}[htb]
    \centering
    \includegraphics[width=0.8\textwidth]{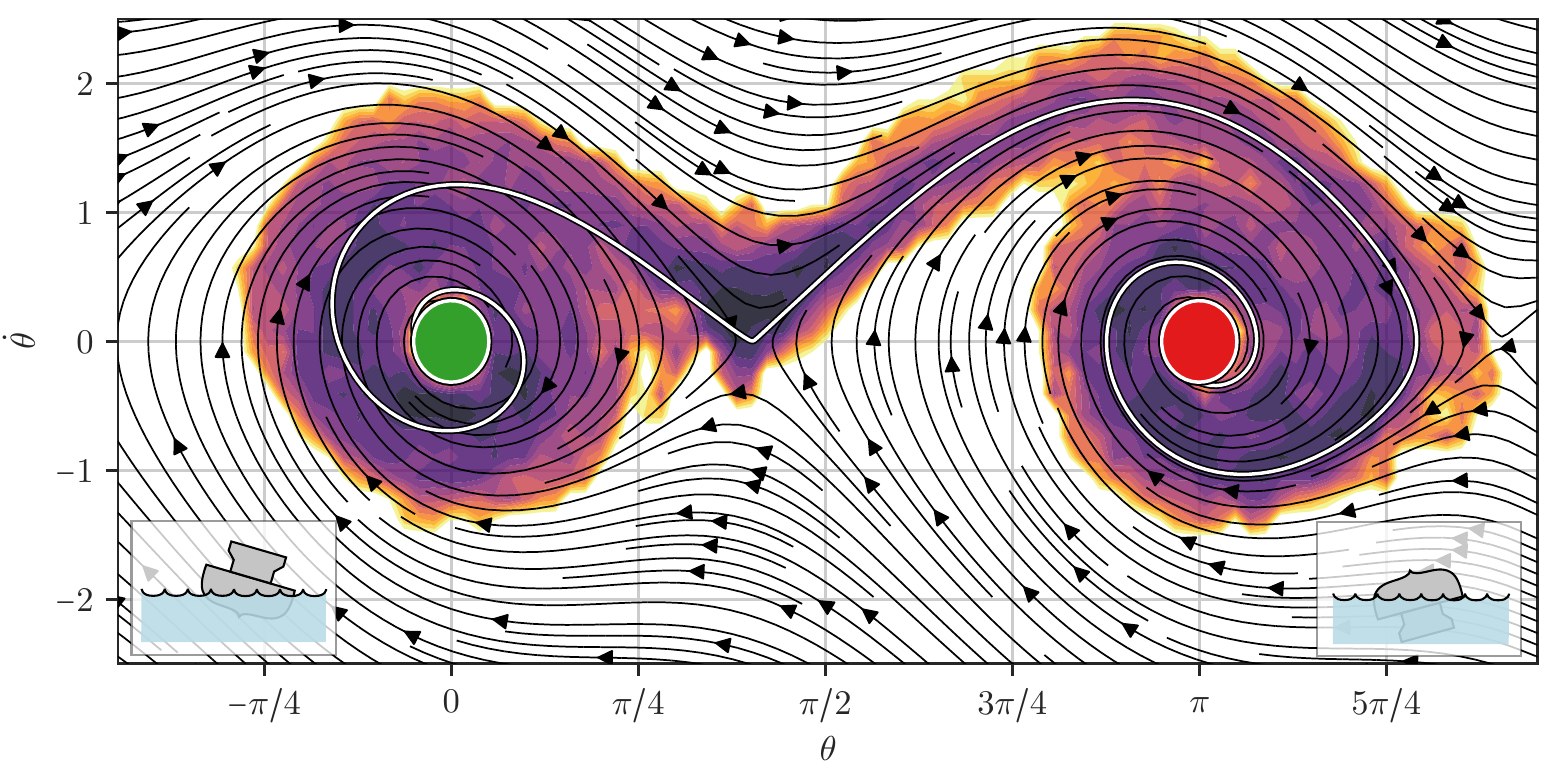}
    \caption{Toy model of a ship capsize as 2-dimensional stochastic system for the roll angle $\theta$ and corresponding angular velocity $\dot\theta$. Trajectories from upright (green) to capsized (red) correspond to transition paths out of the stable basin. The deterministic dynamics are shown as vector field, the density of reactive trajectories as shading, and the large deviation minimiser as white line.}
    \label{fig:LDT}
\end{figure*}

A complementary approach to the above formalism is the perspective of
\emph{stochastic reachability} and \emph{large deviations}. 

Stochastic reachability is a technique used in engineering and computer science to deal with safety issues in a quantitative manner. The objective of stochastic reachability analysis is to evaluate the probabilities associated with dynamic optimization problems. This technique can be used for optimal control under uncertainty, for risk assessment of technical systems, and for safety verification. 
Formally, the system is modelled using a stochastic process (e.g., a Markov chain/process, Wiener process, Gaussian process, or diffusion process) and the unsafe region is modelled as a set in its phase space. Stochastic reachability aims to estimate the probability measure of the set of the trajectories that start in a given set of initial states and reach a target set (a possible unsafe set for the system) in a given time interval. 

There is a close connection to the terminology of chemical reaction kinetics:~a chemical reaction can be viewed a
transition from one locally stable position in state space to
another, driven by the system's stochasticity (for example thermal
noise) and against its typical short-time behaviour. The picture is
that of a random walk in an energy landscape, where a barrier must be
overcome for a reaction to happen. Such transition events are
generally very rare on the timescale intrinsic to the stochasticity,
but waiting long enough one will eventually observe them. A body of
literature is concerned with transition events~\cite{truhlar-garrett-klippenstein:1996}, their
dynamics~\cite{vanden-eijnden:2006} and likelihood~\cite{freidlin-wentzell:2012}. The ultimate question is, of
course, an estimate of the probability of observing a transition, or
equivalently, the \emph{transition rate}.

In the situation of ship capsize, a ship in its upright position can
similarly be considered only locally stable: while a large enough
perturbation will topple it into a capsize, there generally is a
generous region in its 12 dimensional state space where
restoring mechanisms, such as its righting moments, keep it afloat
most of the time. A \emph{transition trajectory} or \emph{reactive
  trajectory} for ship capsize, thus, describes the movement of a ship
in time that, starting in an upright position, will eventually hit an
unsafe region and subsequently capsize, due to a rare influence of its
stochastic components, and generally against its restoring forces. In
this sense, ship capsize can be seen as a first hitting problem, or
\emph{stochastic reachability problem}. Analytical characterizations of the stochastic reachability use equations that link the hitting distribution of the unsafe set with the occupation measure of the safe basin. This is based on the operator methods and Dynkin formula associated to Markov processes. Martingale characterization can be also derived from this equation. 

In general the fact that reactive
trajectories are rare outliers in a usually mechanically stable system
renders their observation by experiment or numerical sampling quite
hard. Crucially, though, their rareness often paradoxically makes them
predictable, which is the core idea behind sample path large
deviations. In this paper, we propose a large deviation methodology to deal with stochastic reachability to provide asymptotic estimates for the probabilities of rare events~\cite{freidlin-wentzell:2012}.

To make the above more concrete, we consider the motion of the ship
$x(t) \in \mathds{R}^{12}$ as introduced in section~\ref{sec:ship-models}
to be a continuous-time Markov process with
stationary distribution $\rho(x)$. Denote by $A$ a neighbourhood of
the upright ship state, and by $B$ the unsafe region ultimately
leading to capsize, for example as specified in
section~\ref{sec:flux-over-saddle}. We can define by $q_+(x)$ the
\emph{forward committor}, i.e.~the probability density over state space
that we will visit $B$ (capsize) before $A$ (righting), or in other
words the probability that we have committed to a capsizing event when
being located at $x$. Similarly, the \emph{backwards committor} is the
probability density that the process at $x$ originates from $A$ rather than
$B$. Given these, the density of reactive trajectories is immediately
available as $\rho_{R} = q_+ \rho q_-$, as can be intuited by reading
the formula as the combined probability of coming from $A$, being at
$x$, and going to $B$. From committor functions and the density of
reactive trajectories, one can finally deduce quantities such as the
probability flux $j_{AB}$ towards capsize, and the capsize rate
$k_{AB}$, with specific formulas depending on the nature of the
process.

The above quantities are generally not accessible for any system of
interest as they necessitate the solution of Dirichlet boundary value problems
(similar to the Fokker-Planck equation) in high
dimensions. Fortunately, this becomes drastically better in the
presence of a \emph{large deviation principle} (LDP). Intuitively
speaking, one obtains the probability of observing an outcome by
integrating (or summing) over all possible ways this outcome can
occur, weighted by its respective probability. The same is true for
reactive trajectories by defining an appropriate path measure. In the
presence of an LDP, this integral can be replaced in an appropriate
limit (such as thermodynamic limit, low temperature limit, or small
noise limit) by the value of the integrant at the most likely path
realizing the outcome. In essence, the integral is computed by a
Laplace method, exchanging a costly transition sampling problem with
an optimisation problem. Knowledge of the large deviation optimal path
allows the computation of transition rates in the large deviation
limit, and the optimal path can be computed quite efficiently by numerical
means even for rather complex systems~\cite{grafke-vanden-eijnden:2019}.

This is exemplified in figure~\ref{fig:LDT}: Here, we consider a toy model for ship capsize for the roll angle and its velocity, $(\theta,\dot\theta)\in\mathds{R}^2$, under Gaussian forcing. We want to consider only trajectories leading to capsize, i.e.~starting upright (green set) and ending capsized (red set). While the direction of the righting moment in phase space is depicted as streamlines, the density of reactive trajectories is shown as shading, and the large deviation minimizing trajectory is highlighted in white. It is clear from the picture that the capsize trajectories concentrate around the optimal path predicted by large deviation theory.

\section{Interconnections}
\label{sec:interconnections}

The two strands are closely related.  They both represent the uncertainties in forcing by probability distributions. 
They formulate capsize as transition across some distinguished set, random (but highly correlated with the forcing function) in the first approach, and fixed or not needing to be specified precisely in the second.  They both aim to produce estimates or bounds on the capsize rate, particularly in the regime where it is expected to be rare.

The two descriptions overlap if the forcing functions are assumed to be the result of filtering a white noise, as is often assumed in the literature~\cite{filter_models:2015}.  This means there is a filter state $z \in \R^k$, some $k \in \N$, satisfying in the simplest case $\dot{z}=Az + \eps \xi$, where $A$ is an asymptotically stable matrix, $\xi$ is a multidimensional white noise (say stationary Gaussian) with autocorrelation $\langle \xi(t)\xi(s)^T\rangle = C\delta(t-s)$ for some positive semi-definite matrix $C$ and a small parameter $\eps$.  Then the ship dynamics can be taken to be of the form $\dot{x}=G(x,z)$, where $x$ represents the 12 dimensions of the ship state-space. The probability distribution on the functions $z$ is easy to handle (linear stochastic process), so one could hope to obtain probabilistic results for the flux over a saddle approach.  Considered as a system on $(x,z)$ the model also fits in the Markovian context of the second approach.  Thus the two can be directly compared.


Our hope is that further understanding will allow development of large deviation theory to more general probability distributions over forcing functions, thereby escaping the Markovian restriction of the second approach.

\section{Conclusion}
\label{sec:conclusion}

We have presented a new stochastic framework for studying ship capsize.  It has two parallel strands, both based on transition state theory, one starting from a deterministic view, the other from a Markovian view.  For filtered white noise models of forcing, the two approaches can in principle be carried to conclusion.  A synthesis is required to treat more general probability distributions for forcing functions.
Extensions are required to pass from probability distributions for sea states to those for forcing functions. Once established, this framework could be used as a building-block for the formulation of safety criteria, optimizing vessel design, and control strategies for the captain to avoid capsize.

\bibliography{bibtex}

\end{document}